\documentclass[12pt]{amsart}
\usepackage{amsmath,amssymb,amsbsy,amsfonts,amsthm,latexsym,
                amsopn,amstext,amsxtra,euscript,amscd}
\usepackage{setspace}
\usepackage{mathptmx}

\renewcommand\bar\overline

\newtheorem{thm}{Theorem} 
\newtheorem{lemma}{Lemma}

\def\average{{\lower .8pt\hbox{$\sim$}\kern -10.5pt\sum}}
\def\ve{\varepsilon}

\def\mod{\text{mod\,}}

\date{}

\begin{document}

\title
{Moebius-Walsh correlation bounds and an estimate of Mauduit and Rivat}
\author
{J.~Bourgain}
\address
{Institute for Advanced Study\\
1 Einstein Drive\\
Princeton, NJ 08540}

\thanks{This research was partially supported by NSF grants DMS-0808042 and DMS-0835373.}
\begin{abstract}
We establish small correlation bounds for the Moebius function and the Walsh system, answering affirmatively a question posed by G.~Kalai [Ka].
The argument is based on generalizing the approach of Mauduit and Rivat [M-R] in order to treat Walsh functions of `large weight', while the `small weight'
case follows from recent work due to B.~Green [Gr].
The conclusion is an estimate uniform over the full Walsh system.
A similar result also holds for the Liouville function.
\end{abstract}
\maketitle

\bigskip
\noindent
{\bf \S0. Introduction}

Fix a large integer $\lambda$ and restrict the Moebius function $\mu$ to the interval $[1, 2^\lambda]\cap \mathbb Z=\Omega$.
Identifying $\Omega$ with the Boolean cube $\{0, 1\}^\lambda$ by binary expansion $x= \sum_{0\leq j<\lambda} x_j 2^j$, the Walsh system
$\big\{w_A; A\subset \{0, \ldots, \lambda-1\}\big\}$ is defined by $w_\phi=1$ and
$$
w_A(x)=\prod_{j\in A}(1-2x_j)= e^{i\pi\sum_{j\in A} x_j}.\eqno{(0.1)}
$$
The Walsh functions on $\Omega$ form an orthonormal basis \big(the character group of $(\mathbb Z/ 2\mathbb Z)^\lambda$\big) and given a function $f$ on $\Omega$, we write
$$
f =\sum_{A\subset \{0, \ldots, \lambda-1\}} \hat f(A) w_A\eqno {(0.2)}
$$
where $\hat f(A) =2^{-\lambda}\sum_{n\in\Omega} f(n) w_A(n)$ are the Fourier-Walsh coefficients of $f$.
Understanding the size and distribution of those coefficients is well-known to be important to various issues, in particular in complexity theory
and computer science.
Roughly speaking, a $F-W$ spectrum which is `spread out' indicates a high level of complexity for the function $f$.
We do not elaborate on this theory here and refer the reader to the extensive literature on the subject; see also the preprint of B.~Green [Gr], 
which motivated this Note.

Returning to the Moebius function and the so-called `Moebius randomness law' it seems therefore reasonable to expect that $\mu|_\Omega$ will have a $F-W$ spectrum
that is not localized.
More precisely, we establish the following uniform bound on its $F-W$ coefficients, answering affirmatively a question posed by G.~Kalai.

\begin{thm}
For $\lambda$ large enough,
$$
\max_{A\subset\{0, \ldots, \lambda-1\}}\Big|\sum_{n<2^\lambda} \mu(n) w_A(n) \Big| < 2^{\lambda-\lambda^{1/10}}\eqno{(0.3)}
$$
\end{thm}
(a similar estimate is also valid for the Liouville function).

The proof of (0.3) involves different arguments, depending on the size $|A|$.
Roughly speaking, one distinguishes between the case $|A|=o(\sqrt\lambda)$ and $|A|\gtrsim \sqrt\lambda$.
In the first case, B.~Green already obtained an estimate of the type (0.3), see [Gr].
Part of the technique used in [Gr] is borrowed from Harman and Katai's work [H-K] on prescribing binary digits of the primes.
Let us point out that in this range the problem of estimating the correlation of $\mu$ with a Walsh function is reduced to estimates on the usual Fourier spectrum
of $\mu$ (by an expansion of $w_A$ in the trigonometric system).
The latter is then achieved either by means of Dirichlet $L$-function theory (when the argument $\alpha$ is close to a rational $\frac aq$ with sufficiently
small denominator $q$) or by Vinogradov's estimate when $q$ is large.
At the other end of the spectrum, when $A=\{0, \ldots, \lambda\}$,  Mauduit and Rivat proved that
$$
\Big|\sum_{n< 2^\lambda} \Lambda(n) \widehat {w_A} (n)\Big| < 2^{(1-\ve)\lambda}
\eqno{(0.4)}
$$
for some $\ve>0$.

Here $\Lambda(n)$ stands for the Van Mangold function ([M-R]).
Their motivation was the solution to a problem of Gelfond on the uniform distribution of the sum of the binary digits of the primes.
Of course, their argument gives a similar bound for the Moebius function as well.
Thus
$$
\Big|\sum_{n< 2^\lambda} \mu(n) \widehat w_{\{0, \ldots, \lambda-1\}}(n)\Big| < 2^{(1-\ve)\lambda}.\eqno{(0.5)}
$$
A remarkable feature of the [M-R] method is that the usual type-I, type-II sum approach in the study of sums
$$
\sum_{n<X} \Lambda(n) f(n) \text { or } \sum_{n<X} \mu(n) f(n)
$$ 
is applied directly to $f =w_{\{1, \ldots, \lambda\}}$ without an initial conversion to additive characters (as done in [H-K] and [Gr]).
The main idea in what follows is to generalize the Mauduit-Rivat argument in order to treat all Walshes $w_A$ provided $A$ is not to small
(the latter case being captured by [Gr]).

Needless to say, the $2^{-\lambda^{1/10}}$-saving in (0.3) can surely be improved (this is an issue concerning the treatment of low-weight Walsh functions)
and no effort has been made in this respect.
We also observe that, assuming $GRH$, (0.3) may be improved to

\begin{thm}
Under $GRH$, assuming $\lambda$ large, we have
$$
\max_{A\subset\{0, \ldots, \lambda-1\}} \Big|\sum_{n< 2^\lambda} \mu(n) w_A(n)\Big|< 2^{\lambda\big(1-\frac c{(\log \lambda)^2}\big)}.
\eqno{(0.6)}
$$
\end{thm}

We will assume the reader familiar with the basic technique, going back to Vinogradov, of type-I and type-II sums, to which sums $\sum_{n<X}\mu(n) f(n)$ may be 
reduced; see [I-K] or [M-R].
In fact, we will rely here on the same version as used in [M-R] (see [M-R], Lemma 1).
Otherwise, besides referring to the work of B.~Green for $|A|$ small, our presentation is basically selfcontained.
In particular, all the required lemmas pertaining to bounds on Fourier coefficients of Walsh functions are proven (they include estimates similar to those
needed in [M-R] and also some additional ones) and are presented in \S1 of the paper.

\vfill\eject

\noindent
{\bf 1. Estimates on Fourier coefficients of Walsh functions}

For $A\subset\{0, \ldots, \lambda-1\}$ and $x=\sum_j x_j 2^j\in [1, 2^\lambda]\cap \mathbb Z$
$$
w_A(x) =\prod_{j\in A}(1-2x_j)= e^{i\pi\sum\limits_{j\in A} x_j} =\prod_{j\in A} h\left(\frac x{2^{j+1}}\right)\eqno{(1.0)}
$$
where $h:\mathbb R\to \{1, -1\}$ is the 1-periodic function
$$
\begin{cases}
h=1  &\text { if } 0\leq x<\frac 12\\
h=-1 &\text { if } \frac 12 \leq x< 1
\end{cases}
$$
For $x\in \mathbb Z$, 
$$
h\left(\frac x{2^{j+1}}\right) =\sum\limits_{|r|< 2^{j+1}} a_{r, j} \  e \left(\frac {rx}{2^{j+1}}\right)\text{
with $\sum |a_r|\lesssim j$}.
$$
It follows that

\begin{lemma}
$w_A(x) =\sum\limits_{k< 2^\lambda} \widehat w_A(k) \, e\left(\frac {kx}{2^\lambda}\right)$ with
$$
\sum|\widehat w_A (k)| <(C\lambda)^{|A|}.\eqno{(1.1)}
$$
\end{lemma}

From the second equality in (1.0), also
$$
\widehat w_A(k) = 2^{-\lambda} \sum _{\{x_j\}} 
e^{i\pi \sum\limits_{j\in A} x_j} \, e^{2\pi i \frac k{2^\lambda} \sum x_j 2^j}
=\prod_{j\not\in A} \left(\frac {1+e(k 2^{j-\lambda})}{2}\right) \prod_{j\in A} \left(\frac {1-e (k
2^{j-\lambda})}2\right)
$$
and
$$
|\widehat w_A(k)|= \prod_{j\not\in A}|\cos \pi k2^{j-\lambda}| \prod_{j\in A} |\sin \pi k2^{j-\lambda}|\eqno{(1.2)}
$$

\begin{lemma}
$$
\Vert\widehat w_A\Vert_\infty \lesssim 2^{-c|A|} \text { for some constant $c>0$}.\eqno{(1.3)}
$$
\end{lemma}

\begin{proof}
Use (1.2).

Taking some $i_0\in A$ and assuming
$$
\left|\sin\pi\frac k{2^{\lambda-i_0}}\right|\approx 1, \text { \rm hence }  \ \left\Vert \frac k{2^{\lambda-i_0}} -\frac
12\right\Vert\approx 0
$$
it follows that either
$$
\left\Vert\frac k{2^{\lambda-i_0-1}} -\frac 14\right\Vert \approx 0
$$
or
$$
\left\Vert \frac k{2^{\lambda-i_0-1}} -\frac 34\right\Vert\approx 0
$$
and in either case
$$
\left|\cos \pi \frac k{2^{\lambda-i_0-1}}\right|, \left|\sin \pi \frac k{2^{\lambda-i_0-1}}\right| \sim \frac 1{\sqrt 2}.
$$
The conclusion follows from (1.2).
\end{proof}

In addition to (1.1), we have the bound

\begin{lemma}
$$
\sum_{k< 2^\lambda} |\widehat w_A(k)|\lesssim 2^{(\frac 12-c)\lambda}.\eqno {(1.4)}
$$
for some constant $c>0$.
\end{lemma}

\begin{proof}

We have to estimate
$$
\sum_{k\in\mathbb Z/ 2^\lambda\mathbb Z} \ \, \prod_{i\leq \lambda} \left|\cos\pi \left(\frac {u_i}2+\frac
k{2^{\lambda-i}}\right)\right|\eqno{(1.5)}
$$
where $u_i=1$ if $i\in A$ and $u_i=0$ if $i\not\in A$.

Perform a shift $k\to k+c2^{\lambda-2} + d2^{\lambda-1}$ with $c, d=0, 1$.

This gives
$$
\sum_{k\in\mathbb Z/2^{\lambda-2}\mathbb Z} \ \, \prod_{2\leq i\leq \lambda} \left|\cos\pi \left(\frac {u_i}2+\frac
k{2^{\lambda-i}}\right)\right|.\eqno{(*)}
$$
with
$$
\begin{aligned}
(*)&=\frac 14\sum_{c, d =0,1} \left|\cos\pi\left(\frac{u_0}2+\frac k{2^\lambda}+\frac c4+\frac d2\right)\right| \
\left|\cos\pi\left(\frac{u_1}2 +\frac k{2^{\lambda-1}} +\frac c2\right)\right|\\[6pt]
&=\frac 14 \sum_{c=0, 1} \left(\left| \cos\pi \left(\frac {u_0}2+\frac k{2^\lambda}+\frac c4\right)\right|+
\left|\sin\pi\left(\frac{u_0}2+\frac k{2^\lambda}+\frac c4\right)\right|\right) 
\left|\cos\pi\left(\frac {u_1} 2+\frac k{2^{\lambda-1}}+\frac c2\right)\right|\\[6pt]
&=\frac 14 \left\{\left(|\cos\phi|+ |\sin\phi|\right).\left|\cos\left(\frac{\pi u_1}2+2\phi\right)\right|+\right.\\[6pt]
& \quad \frac 1{\sqrt 2}(|\cos\phi-\sin\phi\left|+\left|\sin\phi+\cos\phi\right|\right) .\left|\sin\left(
\frac{\pi u_1}2 +2\phi\right)\right|\Big\}
\end{aligned}
\eqno{(1.6)}
$$
where $\phi =\pi\left(\frac {u_0}2+\frac k{2^\lambda}\right)$.
Clearly
$$
\begin{aligned}
(1.6)&\leq \frac 14 \left\{ (1+|\sin 2\phi|)^{\frac 12} \left|{\cos\atop\sin} (2\phi)\right|+(1+|\cos 2\phi|)^{\frac
12}\left|{\sin\atop \cos} (2\phi)\right|\right\}\\[6pt]
&\leq \frac 14 \sqrt {2+\sqrt 2}.
\end{aligned}
$$
Iterating, we obtain the bound
$$
\leq \left (\sqrt{2+\sqrt 2}\right)^{\lambda/2}
$$
and hence (1.4).
\end{proof}

\begin{lemma} Let $r<\lambda$, $a=0, 1, \ldots, 2^r-1$.
Then
$$
\sum_{k\equiv a(\mod 2^r)} |\widehat w_A(k)|\lesssim 2^{(\frac 12-c) (\lambda-r)}.\eqno{(1.7)}
$$
\end{lemma}

\begin{proof}
Writing $k=a+2^rk_1$ with $k_1< 2^{\lambda-r}$,
$$
|\widehat w_A(k)|= \prod_{i<\lambda-r} \left|\cos\pi \left(\frac {u_i}2+\frac a{2^{\lambda-i}}+\frac
{k_1}{2^{\lambda-i-r}}\right)\right|\prod_{i\geq\lambda-r} \left|\cos\pi\left(\frac{u_i}2+\frac
a{2^{\lambda-i}}\right)\right|
$$
$$
\leq \prod_{i<\lambda-r} \left(\left|\cos \pi\left(\frac{u_i} 2+ \frac
{k_1}{2^{\lambda-r-i}}\right)\right|+2^{-\lambda+i+r}\right).\qquad\qquad
\eqno{(1.8)}
$$
For fixed $k_1$, denote
$$
B(k_1) =\left\{i< \lambda-r; \left| \cos\pi\left(\frac{u_i}{2}+\frac {k_1}{2^{\lambda-r-i}}\right)\right|<\left(\frac
1{\sqrt 2}\right)^{\lambda-r-i}\right\}
$$
Hence, if $i\not\in B_{k_1}$
$$
\left|\cos\pi\left(\frac {u_i}2+\frac {k_1}{2^{\lambda-r-i}}\right)\right|+ 2^{-\lambda+r+i} <\left(1+\left(\frac
1{\sqrt 2}\right)^{\lambda-r-i}\right) \left|\cos\pi\left(\frac{u_i}2+\frac{k_1}{2^{\lambda-r-i}}\right)\right|
$$
and if $i\in B_{k_1}$
$$
\left|\cos\pi\left(\frac{u_i}2+\frac{k_1}{2^{\lambda-r-i}}\right)\right|+ 2^{-\lambda+r+i} < \left(\frac 1{\sqrt
2}\right)^{\lambda-r-i} \left(1+ 2\left(\frac 1{\sqrt 2}\right)^{\lambda-r-i}\right) \left|\sin\pi\left(\frac{u_i}
2+\frac {k_1}{2^{\lambda-r-i}}\right)\right|.
$$
Thus certainly
$$
|\widehat w_A(k)|\lesssim \sum_{B\subset\{0, 1, \ldots, \lambda-r-1\}} \left(\frac 1{\sqrt 2}\right)^{\sum\limits_{i\in
B}(\lambda-r-i)}
\prod_{\substack {i\not\in B\\ i<\lambda-r}} \left|\cos\pi\left(\frac{u_i}2+\frac{k_1}{2^{\lambda-r-i}}\right)\right|
\prod_{i\in B} \left|\sin\pi\left(\frac{u_i}{2}+\frac {k_1}{2^{\lambda-r-i}}\right)\right|.\eqno{(1.9)}
$$
Given $B\subset[0, \lambda -r - 1[$, define $B_1\subset [0, \lambda-r-1[$ \  as
$$
B_1=(B\cap[u_i=0])\cup (B^c\cap [u_i=1]).
$$
Hence
$$
(1.9) =\sum\limits_B \left(\frac 1{\sqrt 2}\right) ^{\sum\limits_{i\in B}(\lambda-r-i)} |\widehat w_{B_1} (k_1)|.\eqno{(1.10)}
$$
Summation of (1.10) over $k_1< 2^{\lambda-r}$ and using the bound (1.4) with $\lambda$ replaced by $\lambda-r$ clearly
gives (1.7)
\end{proof}

Next, we also need the following `approximation property' for shifts

\begin{lemma}
Let $A\subset [\lambda-\sigma, \lambda]\cap\mathbb Z$.

Then
$$
\sum_{k< 2^\lambda} |\widehat w_A(k)| < C^{(\log\lambda)^2} (2^\sigma)^{\frac 12-c}.\eqno{(1.11)}
$$
Moreover, there is a bounded function $W_A$ on $[0, \lambda]\cap \mathbb Z$ satisfying $|\widehat W_A |\leq |\widehat w_A|$ and
$$
\left( 2^{-\lambda}\sum_{x< 2^\lambda} |W_A(x) -w_A(x)|^2\right)^{1/2} < 2^{-ct}\leqno{(1.12)}
$$
$$
\widehat W_A(k)=0 \ \text { if } \ |k|> 2^{\sigma+t}\leqno{(1.13)}
$$
Here $t\in\mathbb Z$ is a parameter satisfying $C(\log \lambda)^2< t<\frac 12 (\lambda-\sigma)$.
\end{lemma}
\null\vskip -.8in

\begin{proof}
Writing $k=k_0+2^\sigma k_1$ with $k_0< 2^\sigma, |k_1|< 2^{\lambda-\sigma-1}$ and setting again $u_i =1$ if $i\in A$,
$u_i=0$ if $i\not\in A$, we obtain
$$
|\widehat w_A(k)|=\prod_{i<\lambda-\sigma}\left|\cos\pi \left(\frac{k_0+2^\sigma k_1}{2^{\lambda-i}}\right)\right|.
\prod_{\lambda-\sigma \leq i<\lambda} \left|\cos \pi\left(\frac {u_i}2+ \frac
{k_0}{2^{\lambda-i}}\right)\right|\eqno{(1.14)}
$$
$$
=(1.14). |\widehat w_{A-\lambda+\sigma} (k_0)|.\eqno{(1.15)}
$$
where
$$
A-\lambda+\sigma\subset [0, \sigma]\cap\mathbb Z.
$$
We treat (1.14) as in the proof of Lemma 4, obtaining a
bound
$$
|(1.14)| <\sum\limits_{B\subset \{0, 1, \ldots, \lambda-\sigma-1\}} \left(\frac 1{\sqrt 2}\right)^{\sum\limits_{i\in B}(\lambda-\sigma-i)}
|\widehat w_B(k_1)|.\eqno{(1.16)}
$$

From (1.1), certainly
$$
\sum_{k_1< 2^{\lambda-\sigma}} |\widehat w_B(k_1)|<(C\lambda)^{|B|}\eqno{(1.17)}
$$
and substitution of (1.17) in (1.16) implies by (1.15)
$$
\begin{aligned}
\Vert\widehat w_A\Vert_1 &\leq \Vert\widehat w_{A-\lambda+\sigma}\Vert_1 \, . \, \sum\limits_B\left(\frac 1{\sqrt
2}\right)^{\sum\limits_{i\in B}(\lambda-\sigma-i)} (C\lambda)^{|B|}\\[6pt]
&\overset {\text{Lemma 3}}< (2^\sigma)^{\frac 12-c} \ C^{(\log\lambda)^2}
\end{aligned}
$$
which is (1.11).

Next, let $C(\log \lambda)^2 <\rho<\frac 12(\lambda-\sigma)$ and estimate
$$
\sum_{k_1} \ \sum_{\min B\leq \lambda - \sigma-\rho} \left(\frac 1{\sqrt 2}\right)^{\sum\limits_{i\in B}(\lambda-\sigma-i)}
|\widehat w_B(k_1)|\lesssim 2^{-\rho/4}.\eqno{(1.18)}
$$
If
$$
B\subset [\lambda-\sigma-\rho, \lambda-\sigma]\eqno{(1.19)}
$$
we establish a bound on $\widehat w_B(k_1)$.
Write
$$
\left|\widehat w_B (k_1)\right|= \prod_{i<\lambda-\sigma-\rho} \left|\cos \pi\frac
{k_1}{2^{\lambda-\sigma-i}}\right|.\prod_{\lambda-\sigma-\rho\le i<\lambda-\sigma}\left|\cos\pi
\left(\frac{v_i}2+\frac {k_1}{2^{\lambda-\sigma-i}}\right)\right|
$$
with $v_i=0, 1$ if $i\not\in B$, $i\in B$.
Hence, for $4^\rho<k_1< 2^{\lambda-\sigma-1}$
$$
|\widehat w_B(k_1)|\leq \prod_{\rho<j\leq\lambda-\sigma} \left|\cos\pi \frac {k_1}{2^j}\right|
< k_1^{-c}\eqno{(1.20)}
$$
for some $c< 0$, as we verify by dyadic expansion of $k_1$.

It follows that for $4^\rho \leq K_1 < 2^{\lambda-\sigma}$
$$
\sum_{K_1< |k_1|<  2^{\lambda-\sigma}} \left\{\sum_{B(1.19)}\left(\frac 1{\sqrt 2}\right)^{\sum\limits_{i\in B} (\lambda -\sigma- i)}
|\widehat w_B(k_1)| \right\}^2<
$$
$$
< C \sum_{B(1.19)} \sum_{K_1 < |k_1|<2^{\lambda-\sigma}} \left(\frac 1{\sqrt 2}\right) ^{\sum\limits_{i\in B}(\lambda-\sigma-i)} |\widehat w_B(k_1)|^2\qquad
$$
$$
\overset{(1.20)}<  K_1^{-c} \sum\limits_B \left(\frac 1{\sqrt 2}\right)^{\sum\limits_{i\in B}(\lambda-\sigma-i)}\Vert\widehat w_B\Vert_1\qquad
$$
$$
\overset{(1.17)} < K_1^{-c} C^{(\log\lambda)^2}.\qquad\qquad\qquad
\eqno{(1.21)}
$$
Define $W_A$ as Fourier restriction of $w_A$.
More specifically, let
$$
W_A(x) =\sum \eta (k) \widehat w_A(k) \ e\left(\frac {kx}{2^\lambda}\right)\eqno{(1.22)}
$$
where $\eta:\mathbb R \to[0, 1]$ is trapezoidal with $\eta(z)=1$ for $|z|<K_1 2^\sigma, \eta (z)=0$ for $|z|\geq 2K_12^\sigma$.
Hence $\Vert W_A\Vert_\infty \leq 3$ and $\hat W_A(k)=\hat w_A(k)$ for $|k|\leq K_1 2^\sigma, \hat W_A(k)=0$ for $|k|\geq 2K_1 2^\sigma$.

From the preceding
$$
\Vert \widehat W_A -\widehat w_A\Vert^2_2 \leq \sum_{k_0< 2^\sigma}|\widehat w_{A-\lambda-\sigma}(k_0)|^2 \sum_{K_1\leq |k_1|< 2^{\lambda-\sigma}}
(1.16)^2
$$
$$
\overset {(1.18), (1.21)}< 2^{-\rho/2}+ K_1^{-c} \ C^{(\log\lambda)^2}.\qquad\qquad 
\eqno{(1.23)}
$$
Taking $K_1=2^{t-1}, \rho=\frac {t-1} 2$, Lemma 5 follows.
\end{proof}

The role of $W_A$ is to provide a substitute for $w_A$ with localized Fourier transform.

\begin{lemma}
If $J\subset [1, 2^\lambda[$ \ is an interval, there is a bound
$$
\sum_{k\in J} |\widehat w_A (k)|\lesssim |J|^{\frac 12-c}.\eqno{(1.24)}
$$
\end{lemma}

\begin{proof}

Write
$$
|\widehat w_A(k)|= \prod_{i<\lambda} \left|\cos \pi\left(\frac{u_i}2+\frac k{2^{\lambda -i}}\right)\right|
$$
with $u_i=0 (u_i=1)$ if \ $i\not\in A$ $(i\in A)$.

Assume $2^m\sim |J|< 2^m$. Obviously
$$
\begin{aligned}
|\widehat w_A(k)|\leq \prod_{\lambda-m\leq i<\lambda} \left|\cos\pi\left(\frac {u_i}2+\frac k{2^{\lambda-i}}\right)\right| &=\prod_{0\leq i_1< m}
\left|\cos\pi \left(\frac {u_{i_1+\lambda-m}}2 +\frac k{2^{m-i_1}}\right)\right|\\[6pt]
&= |\widehat w_{A_1}(k)|
\end{aligned}
$$
where
$$
A_1 =\{0\leq i_1<m; \  i_1\in A+ m-\lambda\}.
$$
Hence, since $\hat w_{A_1}$ is $2^m$-periodic
$$
\sum_{k\in J} |\widehat w_A(k)|\leq \sum_{k\in J} |\widehat w_{A_1} (k)|\leq \sum_{k< 2^{m}} |\widehat w_{A_1} (k)|\leq \Vert\hat w_{A_1}\Vert_1< 2^{m(\frac 12-c)}
$$
by Lemma 3.
\end{proof}

\bigskip

\noindent
{\bf 2. Type-II sums}

Let $X=2^\lambda$, 
$S\subset \{0, \ldots, \lambda-1\}$,  $w_S(x) =\prod_{i\in S} (1-2 x_i)$ with  $x=\sum x_i 2^i$.

Specify ranges $M\sim 2^\mu, N\sim 2^\nu$ such that $M\leq N$ and  $M.N\sim X$.

Our goal is to bound bilinear sums of the form $\sum_{\substack {m\sim M\\ n\sim N}}\alpha_m\beta_n w_S(m. n)$, where \hfill\break
$|\alpha_m|, |\beta_n|\leq 1$ are arbitrary
coefficients.

We fix a relatively small dyadic integer $L=2^\rho$ (to be specified).
We assume \hfill\break
$\rho<\frac\mu{100}$, noting that otherwise our final estimate (2.29) is trivial.

Following [M-R], we proceed with the initial reduction of the problem, crucial to our analysis.

Estimate
$$
\left|\sum_{\substack {m\sim M\\ n\sim N}} \alpha_m \beta_n w_S(m.n)\right| \leq \sum_{m\sim M} \left|\sum_{n\sim N} \beta_n w(m.n)\right|.\eqno{(2.1)}
$$
Fix $K$, such that $L2^K< N$ and write using Cauchy's inequality
$$
\begin{aligned}
&\left|\sum_{n\sim N}\beta_n w(m.n)\right| \leq\frac 1L \sum_{n\sim N} \left|\sum^L_{\ell=1} \beta_{n+ \ell2^K}\ w\big(m(n+\ell 2^K)\big)\right|\\[6pt]
&\left|\sum_{n\sim N}\beta_n w(m.n)\right|^2 \lesssim \frac NL \left[\sum_{\substack{n\sim N\\ |\ell|<L}} \beta_n.\overline\beta _{n+\ell.2^K} w(m.n) \ w\big(m(n+\ell
2^K)\big)\right].
\end{aligned}
$$
Hence, by another application of Cauchy's inequality, we obtain
$$
(2.1)^2 \lesssim \frac {M.N} L \sum_{\substack{n\sim N\\ |\ell|<L}} \left|\sum_{m\sim M} w_S(m.n) \ w_S\big(m(n+\ell 2^K)\big)\right|.\eqno{(2.2)}
$$
Comparing the binary expansions of $mn$ and $mn+\ell m 2^K$,
the $K$ first digits remain and  we can assume that also digits $j> K+\mu+\rho+\ve\rho$ are unchanged provided in (2.2) we introduce an additional error 
term of the order $2^{-\ve\rho} M^2N^2 $ (cf. Lemma 5 in [M-R]).
Here $\ve>0$ remains to be specified and we assume $\ve\rho\in\mathbb Z_+$.

Therefore we may write, up to above error
$$
w_S(mn) w_S\big(m(n+\ell 2^K)\big)  \text{ `=' }  w_{S'} (mn) w_{S'} \big(m(n+\ell 2^K)\big)
$$
with
$$
S'=S\cap [K, K+\mu+\rho'] \text { and } \rho'=(1+\ve)\rho
$$
and in (2.2) we may replace $w=w_S$ by $w_{S'}$.

We will either choose $K=0$ or $\mu-\rho\leq K<\lambda-\mu-\rho$.
Hence, by varying $K$, the intervals $[K, K+\mu+\rho]$ will cover $[0, \lambda[$.

For $K\not= 0$, we approximate $w_{S'}$ by $W_{S'}$ given by Lemma 5, applied with $\lambda$ replaced by $K+\mu+\rho'$ and $\sigma$ by $\mu+\rho'$.

Take $t=\ve\rho$ where $\rho$ is certainly assumed to satisfy
$$
\frac \mu{100}> \rho\gg (\log\lambda)^2.
$$
Thus from (1.12)
$$
\sum_{x<X} |w_{S'} (x) -W_{S'} (x)|^2 < 2^{-ct} X.
$$
From the preceding (since $W_{S'}$ is bounded)
$$
(2.2) \lesssim \frac XL \sum_{\substack {n\sim N\\ |\ell|<L}} \left|\sum_{m\sim M} W_{S'}(m.n) W_{S'}\big(m(n+\ell 2^K)\big)\right| \eqno{(2.3)}
$$
$$
+X\sum_{\substack{m\sim M\\ n\sim N}} \left|w_{S'} (mn) -W_{S'}(m.n)\right| \quad \eqno{(2.4)}
$$
$$
+ X^2L^{-\ve}\qquad\qquad\qquad\qquad\qquad
$$
where
$$
\begin{aligned}
(2.4) &< X\left(\sum_{x<X} |w_{S'} (x) -W_{S'}(x)|^2\right)^{\frac 12} \left(\sum_{x\in X} d(x)^2\right)^{\frac 12}\\[6pt]
& <L^{-c\ve} X^2 (\log X)^C < L^{-c\ve} X^2.
\end{aligned}
$$
For $K=0$,
$$
w_{S'}(x) =\sum_{k< 2^{\mu+\rho'}} \widehat w_{S'} (k) \, e\left(\frac {kx}{2^{\mu+\rho'}}\right)\eqno{(2.5)}
$$
where, from Lemma 2 and Lemma 3 applied with $\lambda$ replaced by $\mu+\rho'$
$$
\Vert\widehat w_{S'}\Vert_\infty < 2^{-c|S'|}\eqno{(2.6)}
$$
and
$$
\Vert\widehat w_{S'} \Vert_1 < 2^{(\frac 12 -c)(\mu+\rho')} <2^{(\frac 12-c)(\mu+\rho)}\eqno{(2.7)}
$$
for $\ve$ small enough. 

For $K\not=0$,
$$
W_{S'} (x) =\sum_{|k|< 2^{\mu+\rho'+t}} \widehat W_{S'} (k) \, e\left(\frac {kx}{2^{\mu+\rho' +K}}\right)\eqno{(2.8)}
$$
where
$$
\Vert\widehat W_{S'}\Vert_\infty \leq \Vert\widehat w_{S'}\Vert_\infty < 2^{-c|S'|}\eqno{(2.9)}
$$
and by (1.11) and our choice of $\rho$
$$
\Vert\widehat W_{S'}\Vert_1 < 2^{(\frac 12 -c)(\mu+\rho)}.\eqno{(2.10)}
$$
Denoting by $w$ either $w_{S'}$ when $K=0$ or $W_{S'}$ for $\mu+\rho\leq K<\lambda-\mu-\rho$, substitution of
(2.5), (2.8) and applying a smoothened $m$-summation gives for (2.3), with $M_1=M^{1-\ve_1}$
$$
\frac {M^2 N}L \sum_{\substack {|\ell|\lesssim L\\ n\sim N}} \  \, \sum_{k, k'} |\widehat w(k)| \ |\widehat w(k')| 
1_{\big[\Vert\frac {kn}{2^{\mu+\rho'+K}} -\frac {k'(n+\ell 2^K)} {2^{\mu+\rho'+K}}\Vert < \frac 1{M_1}\big]}\eqno{(2.11)}
$$
up to a negligible error term.

The condition
$$
\left\Vert\frac{(k-k')n}{2^{\mu+\rho'+K}}-\frac {k'\ell}{2^{\mu+\rho'}}\right\Vert < \frac 1{M_1}\eqno{(2.12)}
$$
has to be analyzed.

For $k=k'$ the contribution is
$$
\frac {M^2N^2}L \sum_{|\ell|\lesssim L}  \ \sum_{|k|< 2^{\mu+\rho'+t} } |\widehat w(k)|^2 \ 1_
{\big[\Vert\frac {k\ell}{2^{\mu+\rho'}}\Vert
< \frac 1{M_1}\big]}.\eqno{(2.13)}
$$
The $\ell=0$ contribution in (2.2) is at most  $\frac {M^2N^2}L$.

For $\ell\not= 0$, we get a bound 
$$
M^{2+\ve_1}  N^2  2^{\rho'+t} \Vert \widehat w\Vert_\infty^2 < M^2 N^2 L^2\vert\widehat w\Vert^2_\infty < X^2 L^2 2^{-c|S'|}\eqno{(2.14)}
$$
from (2.6), (2.9) and choosing $\ve_1>0$ small enough to ensure $\ve_1\lambda<\ve\rho$.

In the sequel, we assume $k\not= k'$, $\ell \not=0$.

Also, if in (2.11) for given $k, k',\ell$ there are at most $O(1)$ values of $n$ satisfying (2.12), the resulting contribution is at most
$$
M^2N\Vert\widehat w\Vert^2_1 \operatornamewithlimits<\limits_{\substack {(2.7)\\ (2.10)}} M^2 N(ML)^{1-2c}< X^2 LN^{-c}\eqno{(2.15)}
$$
since $M\leq N$.

Returning to (2.11), consider first the case $K=0$.

We estimate the contribution for
$$
(k-k', 2^{\mu+\rho'})=2^r.
$$
Thus $k-k'=k_1 2^r$, $(k_1, 2)=1$ and (2.12) becomes
$$
\left\Vert\frac {k_1 n}{2^{\mu+\rho'-r}}- \frac {k'\ell}{2^{\mu+\rho'}} \right\Vert < \frac 1{M_1} \eqno{(2.16)}
$$
implying also
$$
\left\Vert \frac{k'\ell}{2^r}\right\Vert < \frac {L^{1+2\ve}}{2^r}.\eqno{(2.17)}
$$
It follows from (2.17) that there are at most $L^{1+2\ve}$ possibilities for $k'$ $(\mod 2^r)$ and hence
for $(k, k')$ $(\mod 2^r)$.

For fixed $k, k', \ell$, (2.16) determines $n$ $(\mod 2^{\mu+\rho'-r})$ up to $1+L^{1+2\ve} 2^{-r}$ possibilities and hence $n$ 
up to $\frac {N2^r}{ML}(1+L^{1+2\ve} 2^{-r})$ possibilities.

Thus the corresponding contribution to (2.11) is at most
$$
\begin{aligned}
&\frac {M^2N}L \sum_{|\ell|\lesssim L} L^{1+2\ve} \frac {N2^r}{ML} (1+L^{1+2\ve} 2^{-r}) \max_{a} \sum_{\substack{k\equiv a (\mod 2^r)\\ k'\equiv a(\mod 2^r)}}
|\widehat w(k)| \ |\widehat w(k')|\\[6pt]
&\lesssim MN^2 (L+2^r) L^{2\ve} \max_a \Bigg[\sum_{\substack {k< 2^{\mu+\rho'}\\ k\equiv a(\mod 2^r)}} |\widehat w(k)| \Bigg]^2.
\end{aligned}
\eqno{(2.18)}
$$
From Lemma 4 applied with $\lambda$ replaced by $\mu+\rho'$
$$
\begin{aligned}
(2.18) &\lesssim  MN^2 (L+2^r)(2^{\mu+\rho' -r})^{1-c} L^{2\ve}\\
&=M^2 N^2 (L^2 2^{-r} +L)(ML2^{-r})^{-c} L^{3\ve}.
\end{aligned}\eqno{(2.19)}
$$
Hence, assuming
$$
ML2^{-r} >L^C\eqno{(2.20)}
$$
we obtain the bound
$$
\frac {X^2}L.
$$
Next, assume
$$
ML2^{-r} < L^C.\eqno{(2.21)}
$$
From the preceding, there are at most $L^{1+4\ve} (ML2^{-r})^2 <L^C$  possibilities for $(k, k')$.

This gives the contribution
$$
M^2N^2 L^C \Vert\widehat w\Vert^2_\infty < L^C X^2 2^{-c|S'|}
$$
and in conclusion $(K=0)$ the bound
$$
X^2(L^{-1} +L^C 2^{-c|S'|}).\eqno{(2.22)}
$$

Next, assume
$$
 K\geq \mu-\rho.\eqno{(2.23)}
$$
Return to (2.11).
Fix $\ell, k, k'$ with $|k-k'|\sim \Delta k<ML^2$.
Letting $n$ range over an interval of size $\frac {ML 2^K}{\Delta k}$, the number of possibilities for $n$ in that interval is at most
$$
1+\frac{L^{1+2\ve} 2^K}{\Delta k}.
$$
Assume
$$
N\gtrsim \frac {ML2^K}{\Delta k}.
$$

The number of $n$'s satisfying (2.12) is at most \big(since $L2^K\geq M>\frac {\Delta K}{L^2}$  by (2.23)\big)
$$
\frac {N\Delta k}{ML2^K} \left(1+\frac {L^{1+2\ve}2^K}{\Delta k}\right)  < \frac NM L^2.
$$
This gives the contribution in (2.11)
$$
L^2 MN^2 \Vert\widehat w\Vert^2_1 \underset{(2.10)}<  L^2MN^2(ML^2)^{1-c} < X^2 L^3M^{-c}.\eqno{(2.24)}
$$

Next, assume
$$
N\ll \frac {ML2^K}{\Delta k}.
$$
From (2.12), for $\ell, k , k'$ given, there are at most
$$
1+\frac{2^KL^3}{\Delta k} \sim \frac {2^KL^3}{\Delta k}
$$
values of $n$.

Also
$$
\left\Vert\frac{k'\ell}{2^{\mu+\rho'}}\right\Vert <\frac 1{M_1}+ \frac {\Delta k.N}{M. 2^{\rho'} .2^K}.
$$
Since $|k'\ell|< 2^{\mu+\rho} L^2$, there is some integer $\ell_1, |\ell_1| < L^2$ s.t.
$$
\left|\frac {k'\ell}{2^{\mu+\rho'}} -\ell_1\right|<\frac 1{M_1}+\frac {\Delta k.N}{M.2^{\rho'} .2^K}
$$
hence
$$
\left|k'-\ell_1 \frac {2^{\mu+\rho'}}{\ell}\right|<L^{1+2\ve}  +\frac {\Delta k.N}{2^K}.
$$
This restricts $k'$ to at most $L^2$ intervals of size $L^{1+2\ve}+\frac {\Delta k.N}{2^K}$.

Using Lemma 6, we obtain the following bound for the contribution to (2.11)
$$
M^2N.L^2 \left(L^{1+2\ve} +\frac {\Delta k.N}{2^K}\right)^{1-c} \frac {2^K{L^3}}{\Delta k}\lesssim
$$
$$
\frac {M^2NL^72^K}{\Delta k}+M^2N^2L^5 \left(\frac {\Delta k.N}{2^K}\right)^{-c} < M^2N^2L^7\left(\frac{2^K}{N.\Delta k}\right)^c.
\eqno{(2.25)}
$$
If we assume
$$
\frac {N.\Delta k}{2^K}>L^C
$$
(2.25) gives the bound
$$
\frac {X^2}L.\eqno{(2.26)}
$$
Assume next
$$
\frac {N.\Delta k}{2^K}<L^C.
$$
From the preceding, $k'$ is restricted to $L^C$ values and the corresponding contribution to (2.11) is bounded by
$$
M^2N^2L^C\Vert\widehat w\Vert_\infty^2 < X^2 L^C 2^{-c|S'|}.\eqno{(2.27)}
$$

Collecting previous bounds gives
$$
(2.11) < X^2 \left(\frac 1L+L^3M^{-c} +L^{C} 2^{-c|S'|}\right)\eqno{(2.28)}
$$
and recalling (2.3), (2.4)
$$
(2.1) <X \left(L^{-c\ve} +L^2M^{-c} +L^C 2^{-c|S'|}\right).\eqno{(2.29)}
$$
In the estimate (2.29), $S'$ depends on the choice of $K$.

Recall that either $K=0$ or $\mu-\rho\leq K<\lambda-\mu-\rho$ and hence, varying $K$, 
the intervals $[K, K+\mu+\rho]$ will cover $[0, \lambda-1]$.
Thus we may choose $K$ as to ensure that
$$
|S'|\geq \max |S\cap J|\gtrsim \frac \mu \lambda |S|\eqno{(2.30)}
$$
with max taken over intervals $J\subset [0, \lambda-1]$ of size $\mu$, in particular (2.29) implies
$$
(2.1) <X\big(L^{-c\ve} +L^2M^{-c} +L^C 2^{-c\frac\mu \lambda|S|}\big)\eqno{(2.31)}
$$
where $L$ is a parameter.

For $|S|\leq \frac {\lambda^{1/2}}H$ with $H\gg 1$ a parameter, we apply B.~Green's estimate (see [Gr])
$$
\left|\sum_{x<2^\lambda} w_S(x) \mu(x)\right|<\lambda e^{-cH}.\eqno{(2.32)}
$$

Thus we assume $|S|>\frac {\lambda^{\frac 12}}H$.
Taking $L=2^H$, it follows from (2.29), (2.31) that
$$
(2.1)\lesssim X.2^{-c\ve H}\eqno{(2.33)}
$$
assuming either that
$$
M>2^{CH^2\lambda^{1/2}}\eqno {(2.34)}
$$
or
$$
M>C^H \text { and } |S'| >CH \ \big(S' \text { satisfying (2.30)}\big).\eqno{(2.35)}
$$

\bigskip

\noindent
{\bf 3. Type-I sums and conclusion}

We use Lemma 1 from [M-R] but treat also some of the type-I sums as type-II sums.
Indeed, according to (2.33), (2.34), only the range $M<C^{H^2\lambda^{1/2}}$ remains to be treated.

Thus we need to bound
$$
\sum_{m\sim M} \left|\sum_{n\sim N} w_S(mn)\right|\eqno{(3.1)}
$$
where $M.N \sim X=2^\lambda, M< C^{H^2\lambda^{1/2}}$.
We assume $|S|>\frac {\lambda^{1/2}}H$.

Expanding in Fourier and using a suitable mollifier in the $n$-summation, we obtain
$$
(3.1)\leq \sum_{m\sim M} \ \sum_{k<X} |\widehat w_S(k)| \ \left|\sum_{n\sim N} \, e\left(\frac {kmn}{2^\lambda}\right)\right|
$$
$$
<N \sum_{\substack {m\sim M\\ k<X}} |\widehat w_S(k)| \,  1_{\big[\Vert\frac{km}{2^\lambda}\Vert < \frac {\lambda^2}N\big]} 
+o(1) \quad\eqno{(3.2)}
$$
$$
< NM^2 \lambda^2\Vert\widehat w_S\Vert_\infty \ \qquad\qquad\qquad \eqno{(3.2')}
$$
$$
 < XM 2^{-c\lambda^{1/2}H^{-1}}\lambda^2.\qquad\qquad \ \eqno{(3.3)}
$$
Taking $H<\lambda^{1/10}$, (3.3) is certainly conclusive if $M<C^H$. Hence
recalling (2.35), we can assume that
$$
\mu> H \text { and } \ \max |S\cap J|<CH\eqno{(3.4)}
$$
for any interval $J\subset \{0, \ldots, \lambda-1\}$ of size $\mu$, where $M\sim 2^\mu$.

Assumption (3.4) will provide further information on $\hat w_S$ that will be useful in exploiting (3.2).

Write 
$$
S=S_1\cup S_2
$$
where $S_1 =S\cap [0, \lambda-2\mu]$ and $S_2 =S\cap [\lambda-2\mu, \lambda]$.
Hence by (3.4),
$$
|S_2|<CH.
$$
Thus
$$
w_{S_2}(x) \overset{(1.0)}= \prod_{j\in S_2} h\left(\frac x{2^{j+1}}\right)
$$
$$
\qquad = \sum_{k_2\in \mathcal A_2} \widehat w_{S_2} (k_2) \, e\left(\frac{k_2x}{2^\lambda}\right)+ O_{L^1} (2^{-H})
\eqno{(3.5)}
$$
where the set $\mathcal A_2$ may be taken of size
$$
|\mathcal A_2| < 2^{H|S_2|} <C^{H^2}\eqno{(3.6)}
$$
(obtained by truncation of the Fourier expansion of $h$).

On the other hand
$$
w_{S_1}(x) =\sum_{k_1< 2^{\lambda-2\mu}} \widehat w_{S_1} (k_1) \, e\left(\frac {k_1 x}{2^{\lambda -2\mu}}\right)
$$
and hence
$$
w_S(x) =\sum_{\substack {k_1< 2^{\lambda-2\mu}\\ k_2\in\mathcal A_2}} \widehat w_{S_1} (k_1) \widehat w_{S_2}(k_2) \,
e\left(\frac {2^{2\mu }k_1+k_2}{2^\lambda} x\right) + O_{L^1}(2^{-H}).\eqno{(3.7)}
$$

The bound (3.2) becomes now
$$
N\sum_{\substack {m\sim M\\ k_1< 2^{\lambda-2\mu}\\ k_2\in \mathcal A_2}} |\widehat w_{S_1} (k_1)| \ |\widehat w_{S_2} (k_2)| 
\ 1_{\big[\big\Vert\frac {2^{2\mu}k_1+k_2}{2^\lambda} m\big\Vert<\frac {\lambda^2}N\big]}
$$
$$
< N|\mathcal A_2| \ \Vert\widehat w_{S_1}\Vert_\infty . \max_{k_2} \sum_{m\sim M} \left| \left\{ k_1< 2^{\lambda-2\mu}; \left\Vert
{\frac {2^{2\mu}k_1+k_2}{2^\lambda}} m\right\Vert <\frac {\lambda^2}N\right\}\right|.\eqno{(3.8)}
$$
Clearly
$$
\begin{aligned}
&\sum_{m\sim M} \left|\left\{ k_1< 2^{\lambda-2\mu}; \left\Vert\frac{k_1m}{2^{\lambda-2\mu}}\right\Vert <\frac {2\lambda^2}
N\right\}\right|=\\[8pt]
&\sum_{m\sim M} \left|\left\{ k_1<2^{\lambda-2\mu}; k_1 m\equiv 0 (\mod 2^{\lambda-2\mu})\right\}\right|\lesssim \mu. M
\end{aligned}
$$
and therefore, since $|S_1| \gtrsim \frac {\lambda^{1/2}}H $ and (3.6)
$$
(3.8) <\mu C^{H^2} 2^{-c\lambda^{1/2} H^{-1}} NM 
$$
$$
\, < 2^{-c\lambda^{1/2} H^{-1}} X.\qquad  \eqno{(3.9)}
$$
From (2.33) and (3.9), we can claim a uniform bound
$$
\left|\sum_{x<X} \mu(x) w_S(x)\right| \lesssim X. 2^{-c\lambda^{1/10}}\eqno{(3.10)}
$$
hence obtaining Theorem 1.
\bigskip

Under GRH, (3.10) can be improved of course.

First, from a result due to Baker and Harman [B-H], there is a uniform bound
$$
\left\Vert\sum_{n\in X} \mu(n) e(n\theta)\right\Vert_\infty \ll X^{\frac 34 +\ve}.\eqno{(3.11)}
$$
Hence
$$
\left|\sum_{n<X}\mu(n) w_S(n)\right|< \Vert\hat w_S\Vert_1\,  X^{\frac 34+\ve'}< (\log X)^{|S|}\,  X^{\frac 34+\ve'}\eqno{(3.12)}
$$
and we may assume
$$
|S|>c \, \frac{\log X}{\log\log X}.\eqno{(3.13)}
$$
If (3.13), apply the type-I-II analysis above.

From (2.31), assuming
$$
M\sim 2^\mu>X^{c_1 \frac 1{\log\log X}}\eqno{(3.14)}
$$
and choosing $L$ appropriately, we obtain
$$
(2.1)< X. 2^{-c\frac{\log X}{(\log\log X)^2}}.\eqno{(3.15)}
$$
If $M$ fails (3.14) the type-I bound (3.2') gives
$$
\begin{aligned}
(3.1)&< X.M \Vert\hat w_S\Vert_\infty\\[6pt]
&\overset{(1.3)}< X.X^{c_1\frac 1{\log\log X}} \ 2^{-c' \frac {\log X}{\log\log X}}\\
\end{aligned}
$$
$$
< X^{1-c_2 \frac 1{\log\log X}}
\qquad
\eqno{(3.16)}
$$
for appropriate choice of $c_1$ in (3.14).

In either case
$$
\left|\sum_{n<X} \mu(n) w_S(n)\right|< X^{1-\frac c{(\log\log X)^2}}\eqno{(3.17)}
$$
which is Theorem 2.

\end{document}